\documentclass[12pt]{article}

\usepackage[margin=1in]{geometry}
\usepackage{charter}
\usepackage[noblocks]{authblk}

\usepackage[usenames,dvipsnames]{xcolor}
\usepackage{graphicx, pgfplots, subfigure, standalone, pgf, tikz}
\usepackage[inline]{enumitem}
\usepackage[justification=centering]{caption}

\usepackage[sort&compress,numbers]{natbib}

\usepackage{comment, url,float, setspace, xspace}
\pgfplotsset{compat=1.18} 

\usepackage{amsmath, amssymb, amsthm, amsfonts, mathtools, array, bm, xfrac}

\setlist[enumerate]{label=(\roman*.)}

\usepackage[ruled,vlined,linesnumbered]{algorithm2e}

\usepackage[short,c2,nocomma]{optidef}

\usepackage[colorlinks,
            linkcolor=BlueViolet,
            citecolor=blue,
            urlcolor=blue]{hyperref}

\usepackage[nameinlink]{cleveref}

\usepackage{csvsimple, booktabs, longtable, multirow, multicol, adjustbox}	%

\usepackage[xindy, acronyms, nowarn]{glossaries}   %
\makeglossaries

\definecolor{lime}{HTML}{A6CE39}
\DeclareRobustCommand{\orcidicon}{
    \hspace{-3mm}
	\begin{tikzpicture} 
    \draw[lime, fill=lime] (0,0) circle [radius=0.15] node[white] { 
        {\fontfamily{qag}\selectfont \tiny ID} 
    };
	\end{tikzpicture} 
    \hspace{-2mm}
}

\foreach \x in {A, ..., Z}{\expandafter\xdef\csname orcid\x\endcsname{\noexpand\href{https://orcid.org/\csname orcidauthor\x\endcsname}{\noexpand\orcidicon} } }
\newtheorem{definition}{Definition}[section]
\newtheorem{example}{Example}[section]

\newtheorem{remark}{Remark}[section]

\crefname{line}{Line}{Lines}
\crefname{lemma}{Lemma}{Lemmata}
\crefname{theorem}{Theorem}{Theorems}
\crefname{proposition}{Proposition}{Propositions}
\crefname{algorithm}{Algorithm}{Algorithms}
\crefname{equation}{}{}
\crefname{definition}{Definition}{Definition}
\crefname{claim}{Claim}{Claim}
\crefname{corollary}{Corollary}{Corollaries}
\crefname{remark}{Remark}{Remarks}
\crefname{example}{Example}{Examples}
\crefname{figure}{Figure}{Figures}
\crefname{section}{Section}{Sections}
\crefname{table}{Table}{Tables}

\newacronym{LO}{LO}{Linear Optimization}
\newacronym{QO}{QO}{Quadratic Optimization}
\newacronym{MIQO}{MIQO}{Mixed-Integer Quadratic Optimization}
\newacronym{MIO}{MIP}{Mixed-Integer Programming}
\newacronym{MILO}{MILP}{Mixed-Integer Linear Programming}
\newacronym{MINLO}{MINLP}{Mixed-Integer Nonlinear Programming}
\newacronym{AGT}{AGT}{Algorithmic Game Theory}
\newacronym{POA}{PoA}{Price of Aggression}
\newacronym{NE}{NE}{Nash Equilibrium}
\newacronym{POS}{PoS}{Price of Security}
\newacronym{IPG}{IPG}{Integer Programming Game}
\newacronym{DiD}{DiD}{Defence-in-Depth}
\newacronym{CSPM}{CSPM}{Cloud Security Posture Management}
\newacronym{CWPP}{CWPP}{Cloud Workload Protection Platform}
\newacronym{CNP}{CNP}{Critical Node Problem}
\newacronym{CNG}{CNG}{Critical Node Game}

\newcommand{\eg}{\textit{e.g.}}
\newcommand{\ie}{\textit{i.e.}}

\title{ The Critical Node Game\thanks{Most of this work has been conducted while the first and third authors were members of the Canada Excellence Research Chair in ``Data Science for Real-time Decision-making" at Polytechnique Montr\'eal. The CERC support and that of Ericsson Canada Inc. are warmly acknowledged.}  }

\author[1]{Gabriele Dragotto\orcidA{}}
\author[2]{Amine Boukhtouta\orcidB{}}
\author[3]{Andrea Lodi\orcidC{}}
\author[4]{Mehdi Taobane\orcidD{}}

\affil[1]{Department of Operations Research and Financial Engineering, Center for Statistics and Machine Learning, Princeton University, U.S.A.}
\affil[2]{Ericsson Canada Inc., Canada}
\affil[3]{Jacobs Technion-Cornell Institute, Cornell Tech and Technion - IIT, U.S.A.}
\affil[4]{CERC, Polytechnique Montr\'eal, Canada}

\begin{document}
\maketitle

\begin{abstract}
	In this work, we introduce a game-theoretic model that assesses the cyber-security risk of cloud networks and informs security experts on the optimal security strategies.
	Our approach combines game theory, combinatorial optimization, and cyber-security and aims to minimize the unexpected network disruptions caused by malicious cyber-attacks under uncertainty. Methodologically, we introduce the \emph{critical node game}, a simultaneous and non-cooperative attacker-defender game where each player solves a combinatorial optimization problem parametrized in the variables of the other player. Each player simultaneously commits to a defensive (or attacking) strategy with limited knowledge about the choices of their adversary. We provide a realistic model for the critical node game and propose an algorithm to compute its stable solutions, \ie, its Nash equilibria. Practically, our approach enables security experts to assess the security posture of the cloud network and dynamically adapt the level of cyber-protection deployed on the network.
	We provide a detailed analysis of a real-world cloud network and demonstrate the efficacy of our approach through extensive computational tests.
\end{abstract}

\section{Introduction}
Cloud networks are the backbone of the modern internet infrastructure, as they provide the resources, applications, and web services we use daily. In a cloud network, network operators provide ready-to-use resources (\eg, servers and virtual machines) to some end-users; for instance, Dropbox provides cloud storage, and Amazon Web Services provides virtual computing infrastructure. The network operators deploy their physical resources (\ie, servers) in several points of presence (\ie, data centers), and they manage the network on behalf of their customers.
Cloud networks provide businesses and individuals with flexible and scalable access to computing resources and data storage, making managing and sharing information easier. However, because these networks rely on the Internet, remote servers, and third-party software, they are also vulnerable to cyber-attacks and expose organizations to greater external risks.

In this work, we focus on responding to uncertain cyber-attacks, network intrusions performed by a malicious attacker whose strategy is uncertain (\ie, any node in a network can be subject to an attack). We introduce the \gls{CNG}, a \emph{non-cooperative} and \emph{simultaneous} 2-player game between a \emph{defender} (\eg, the network operator) and a malicious \emph{attacker}.
The \gls{CNG} falls into the family of integer programming games, a large class of simultaneous non-cooperative games with complete information where each player solves a parameterized (in their opponents' variables) integer program \citep{koppe_rational_2011,IPGs_2023_Tutorial}. Our game-theoretic framework assesses cloud networks' security posture (\ie, readiness to handle attacks) and provides network operators with an \emph{a priori} security recommendations.
In the \gls{CNG}, the decisions of each player are made by solving a \gls{CNP} \citep{lalou_critical_2018}, \ie, the problem of deciding which nodes to remove (attack or defend) in a given graph (representing the network topology) under some budget constraints.
Each of the two players simultaneously and non-cooperatively decides their strategy by solving an integer program parameterized in the variables of the other player.
While the defender decides which nodes to defend (\eg, which nodes to protect with an extra security level), the attacker decides which nodes to attack. Each player's payoff (\ie, objective function) is parameterized in their opponent's moves and measures the effectiveness of the associated defensive or attacking strategy.

On the one hand, the parameters of the defender's optimization problems model the cloud network topology and implicitly model the importance of each node in the network.
On the other hand, as the attack's dynamics are uncertain, the parameters of the attacker's optimization problem model the \emph{class} of attacks that the defender deems reasonable for a potential attacker. In this sense, the uncertainty derives from the attacker optimization problem representing a range of possible attacks rather than a single attack or a predetermined class of attacks. This approach is, in its motivation, close to the idea of robust optimization \citep{bertsimas_robust_2022}: instead of considering a single attack, we consider the whole set of attacks that could happen, given some minimal knowledge about the attacker's capabilities (\eg, the type of vulnerability exploited). We employ the \gls{NE} \citep{nash_equilibrium_1950} as the standard solution concept for the \gls{CNG}. Nash equilibria are \emph{stable} solutions because no player can deviate without diminishing their payoff. Specifically, we aim to select the \emph{best} \gls{NE}, \ie, the \gls{NE} that maximizes the defender's payoff. In other words, we aim to select the strategy that benefits the defender the most while forcing the attacker to behave as the defender expects.

\begin{figure}[!htpb]
	\centering
	\includegraphics[width=0.6\linewidth]{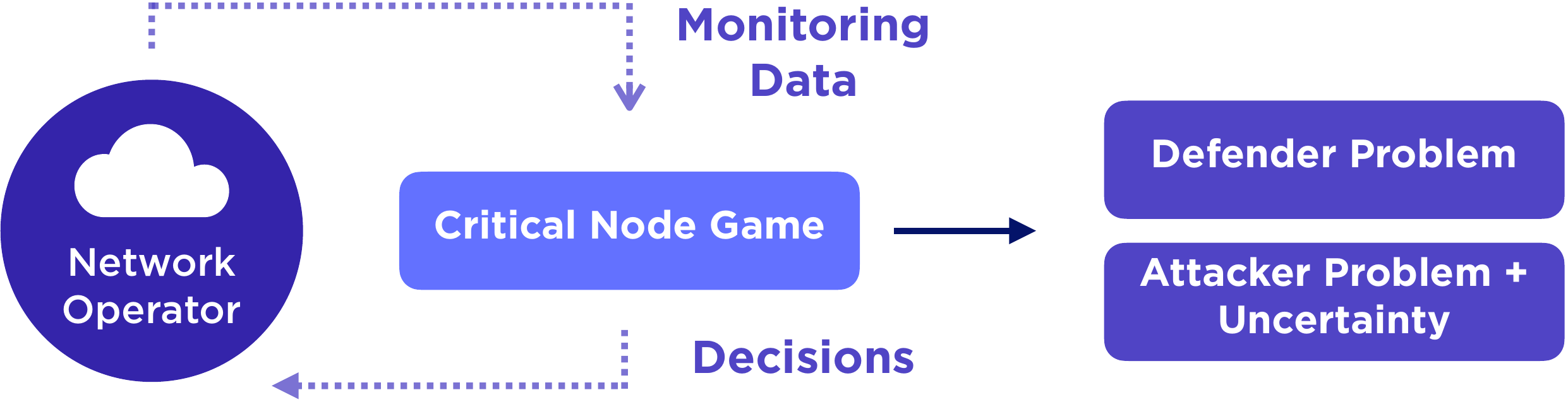}
	\caption{The network operator collects monitoring data that feed the \gls{CNG} model. The \gls{CNG} provides the prescriptive cyber-security recommendations.}
	\label{fig:cycle}
\end{figure}

\paragraph{Our Motivation.}
The cyber-security of cloud networks mainly relies on the so-called \emph{defense-in-depth} approach, \ie, a security paradigm that combines a series of security mechanisms with security best practices to protect the confidentiality, integrity, and availability of the network resources. The core intuition behind defense-in-depth is that a \emph{single} security mechanism is significantly less effective against cyber threats than a stack of \emph{multiple}, diverse and redundant security mechanisms.
This paradigm is implemented throughout cloud networks via several security tools that distribute the security responsibility between the cloud network operators and the end-users. For instance, network operators often maintain network-level security, such as firewalls protecting the whole network and strict user privileges (\eg, the so-called trust policies). In contrast, end-users maintain application-level security, such as custom-tailored software protecting their applications.
The ensemble of tools adapted by operators and users often produces a large amount of \emph{monitoring data} that is generally fed to some \emph{intrusion detection system}. Once this system detects a potential intrusion from a malicious attacker, it generally warns the network operators or puts the system under preventive mitigation (\eg, it increases firewall aggressiveness).
On the one hand, the monitoring data helps the intrusion detection system to identify potential cyber-security threats and security gaps in the cloud infrastructure. On the other hand, however, the detection system is often unable to come up with security prescriptions that emerge from the monitoring data, and it often requires some input from the network operators \citep{roy_survey_2010}.
This is precisely the core motivation behind this work: providing a mathematical framework that uses the monitoring data to quantitatively assess the potential risks associated with cyber-attacks, and dynamically adapt the level of cyber-protection deployed on the network. We represent this information cycle in \cref{fig:cycle}. Finally, we consider a one-shot game instead of a sequential or repeated game by implicitly assuming that both defender and attacker are unaware of each others' strategies. This structural assumption enables the defender to derive prescriptive insights on the defensive strategies with minimal temporal information (\eg, previous attacks) and, thus, more unpredictability on the attacker's strategy (\eg, novel attack techniques).

\paragraph{Our Contribution.}
We summarize our contributions as follows:

\begin{enumerate}
	\item Motivated by the need to assess and react to cyberattacks in highly-complex cloud networks, we propose the \gls{CNG}, a non-cooperative simultaneous 2-player game between a defender and an attacker.  Our model accounts for the uncertainty of the attack by modeling an attacker as an independent decision-maker interacting with the defender. We employ the monitoring data to model the \emph{class of attacks} that are likely to happen instead of deterministically determining which one is more likely to happen. In this sense, our model requires almost no probabilistic assumption to work. In practice, our model enables the network operators to assess the potential risks associated with cyber-attacks and allows the network operator to dynamically adapt the level of cyber protection deployed in the case of an attack.
	\item We formulate the \gls{CNG} as an integer programming game where each player solves a \gls{CNP} formulated as a parametrized binary knapsack problem. We propose two metrics, the \gls{POA} and \gls{POS}, to evaluate the effectiveness of the attacker's and defender's strategies and measure the efficiency of the resulting cyber-security policies. We then tailor the \emph{ZERO Regrets} algorithm of \citet{Dragotto_2021_ZERORegrets} to compute and optimize over the Nash equilibria of the \gls{CNG}.
	\item We provide an extensive set of computational tests on synthetic and real-world instances to demonstrate the capabilities of the \gls{CNG}. Finally, we derive some informative insights from the \gls{CNG} equilibria and demonstrate how they can be practically employed to guide the countermeasures against cyber-attacks.
\end{enumerate}

\paragraph{Outline.} We present a brief literature review in \cref{sec:rw}, while in \cref{sec:model}, we introduce the mathematical formulation, the \gls{CNG}, and the metrics to measure the effectiveness of the defending strategies. In \cref{sec:cng}, we introduce the cloud-based formulation for the \gls{CNG} and briefly describe how to compute its equilibria. In \cref{sec:tests}, we present our experimental setup and the computational tests. Finally, we present our conclusions in \cref{sec:conclusion}.

\section{Literature Review}
\label{sec:rw}

\paragraph{Critical Nodes in Networks.} The idea of detecting the most critical nodes in a network dates back several decades. \citet{ball_finding_1989} studied the problem of removing the most vital arcs in a network to maximize the shortest path among pairs of nodes. Similarly, \citet{bazgan_most_2011} studied the problem of finding the most vital node in a network to maximize the weight of an independent set. \citet{assimakopoulos_network_1987} applied a similar analysis to model an interdiction problem to control the spread of infectious diseases in hospitals. Similar ideas were pioneered by \citet{tao_epidemic_2006,cohen_efficient_2003} in the context of epidemics and computer networks. \citet{commander_wireless_2007} proposed an optimization problem to determine the optimal location of network jamming devices to maximize the network disruption.
\citet{finbow_fireghter_2009} studied the so-called firefighter problem, \ie, the problem of deciding which node to defend in a graph considering a temporal dynamic that spreads a virus or a fire over the nodes.
Several authors investigated the problem of determining the critical node in more classic combinatorial optimization contexts, such as matching graphs, network flows graphs, and in general, in graphs with special structures
\citep{addis_identifying_2013,zenklusen_connectivity_2014,zenklusen_network_2010,zenklusen_matching_2010,di_summa_complexity_2011,arulselvan_detecting_2009,borgatti_identifying_2003,ratliff_finding_1975,shen_polynomial-time_2012,hosteins_stochastic_2020}. We refer to \citet{lalou_critical_2018,zenklusen_matching_2010,veremyev_integer_2014} (and the references therein) for a complete survey on the \gls{CNP}.
Our approach complements the previous literature by extending the \gls{CNP} to a multi-agent non-cooperative setting, as opposed to a single-agent model, and by contextualizing the model in cloud networks cyber-security.

\paragraph{Interdiction Games.} A different yet connected stream of literature considers so-called interdiction games. An interdiction game is a sequential game (\eg, a Stackelberg game) where the first-movers, or the players from the previous rounds, can restrain the next-round players from using some resources. \citet{fischetti_interdiction_2019, taninmis_branch-and-cut_2022} proposed two algorithms to compute the solutions to knapsack interdiction games (\ie, interdiction games where players solve knapsack problems with interdiction constraints) and submodular interdiction games, respectively. \citet{israeli_shortest-path_2002} modeled an interdiction shortest path problem, where interdiction occurs on the edges.
Recently, \citet{baggio_multilevel_2021} proposed a multi-level interdiction approach to the critical node problem, where a set of players sequentially plan infrastructure design over a graph: the approach involves a mathematical program with nested optimization problems and resembles the work of \citet{brown_defending_2006} in the context of critical infrastructures.
\citep{nabli_curriculum_2020} developed a general greedy heuristic for network defender-attacker-defender problems with knapsack constraints using reinforcement learning and validated the approach on the game of \citep{baggio_multilevel_2021}.
\citet{shen_exact_2012} studied the properties of the optimization formulations associated with interdiction problems over undirected graphs.
While the interdiction games literature often assumes players act in rounds, our work approaches the problem from a simultaneous perspective; specifically, we let the attacker and defender play simultaneously without knowing the opponents' strategy. In this sense, our approach is closer to the realm of robust optimization \citep{bertsimas_robust_2022} than to those of interdiction games. Furthermore, compared to \citep{baggio_multilevel_2021,brown_defending_2006,nabli_curriculum_2020}, our \gls{CNG} does not assume that the game is zero-sum, \ie, that the attacker's payoff is the negative defender's payoff.

\paragraph{Game Theory and Cyber-security.} \citet{roy_survey_2010} provided a detailed survey on game-theoretic models for cyber-security. The authors argue that game-theoretical frameworks for network security primarily possess two strengths. First, they can model complex strategic decision-making settings and compare several scenarios before deciding. Second, compared to most of the state-of-the-art heuristic approaches, game theory models are quantitative and exact methods. \citet{he2008game,chen2009game} modeled a similar attacker-defender game to derive insights into the defender's strategy and the resources needed to defend. However, in contrast with those approaches, (i.) we formulate a problem starting from a graph, thus incorporating the structure of the network inside our model, (ii.) we explicitly represent each player's problem as an optimization problem, in the spirit of \citep{Dragotto_2022_Thesis}, and (iii.) we compute and select a specific \gls{NE}. Finally, \citet{goyal_attack_2014} and \citet{dziubinski_network_2013} develop a sequential (\eg, Stackelberg) model that shares a similar payoff structure with our \gls{CNG}, and provide several structural results connecting the graph topology and the defensive strategies.
In contrast to our model, the previous works assume a sequential structure where the network designer plays first, and the attacker follows, whereas we assume players act simultaneously. In game-theoretic terms, this implies that the defender (and the attacker) is unaware of the realization of their opponent's strategy, thus inducing a different information set for each player. The models in \citep{goyal_attack_2014,dziubinski_network_2013} are also closely related to the multi-level plan-attack-defend framework of \citet{baggio_multilevel_2021}. However, the latter considers an interdictive setting.

\section{Mathematical Formulation}
\label{sec:model}

Let $G=(V,E)$ be an undirected graph representing the cloud network topology, where the resources $V$ (\eg, routers, servers, cloud instances) are linked through some connections $E$. Each edge $e \in E$ is made of a pair of nodes $i,j \in V$ so that $(i,j)$ is a connection between the two resources $i$ and $j$. In \cref{def:CNP}, we give an abstract definition of the \gls{CNP}.
\begin{definition}[Critical Node Problem]
	Given a graph $G=(V, E)$ and $k \in \mathbb{N}$, the \gls{CNP} is the problem of removing a subset of nodes $S \subseteq V$ with $|S|<k$ so that $G=(V\backslash S,E)$ minimizes a function of the connectivity of $G$.
	\label{def:CNP}
\end{definition}
In other words, the \gls{CNP} asks to detect the most critical $k$ nodes of $G$ with respect to a given connectivity function. This function can be, for instance, the number of pairwise connected nodes or the number of connected components in $G$ (see \citep{veremyev_integer_2014} for a detailed review of families of functions). Similarly, when $G$ represents a cloud network, the problem equivalently asks to determine the $k$ most critical devices in the network so that removing them, \ie, removing $S \subseteq V$, maximizes the service disruption (\ie, the number of unavailable resources) on $G$. The \gls{CNP} possesses an intuitive formulation as a combinatorial optimization problem, and it is often solved as an equivalent integer program.

We focus on a  \gls{CNP} formulated as a knapsack problem, and we take a dual interpretation, where the selection of the critical nodes is a way of protecting the network, \ie, maximizing the network's available connectivity (in case of an attack).
Let $x_i$ be a variable associated with each node $i \in V$, and let $x_i=1$ if and only if $i$ is protected. The subset $S$ is then the set $\{ i \in V  : x_i=1 \}$. Assuming that a function $f(x)$ measures the connectivity of $S$ induced by $x$, the \gls{CNP} corresponds to the problem
\begin{align}
	\max_x \Big \{ f(x) : w^\top x \le W, x \in \{0,1\}^{|V|} \Big\},
	\label{eq:CNP}
\end{align}
where $w \in \mathbb{R}^{|V|}_+$ and $W \in \mathbb{R}_+$ are the parameters of the weighted knapsack constraint. Whenever $w$ is a vector of ones, the constraint becomes a knapsack constraint asking to remove at most $\lfloor W \rfloor$ nodes from $G$.

\subsection{The Critical Node Game}
While the \gls{CNP} models a wide range of applications from interdiction to network resilience, it only encompasses a single decision-maker. Therefore, we extend the \gls{CNP} of \cref{eq:CNP} into the \gls{CNG} by introducing two different decision-makers: the defender, who controls the variables $x$, and the attacker, who controls the variables $\alpha$. Similarly to the \gls{CNP}, we assume that for any $i \in V$, $x_i=1$ if and only if the defender \emph{protects} node $i$, and $\alpha_i=1$ if and only if the attacker \emph{attacks} node $i$.
\begin{definition}[Critical Node Game]
	Given a graph $G$, the \gls{CNG} is a 2-player \emph{simultaneous} and \emph{non-cooperative} game with \emph{complete information} where the first player (the defender) solves
	\begin{align}
		\max_{x} \Big \{ f^d(x;\alpha) : d^\top x \le D, x \in \{0,1\}^{|V|} \Big\},
	\end{align}
	and the second player (the attacker) solves
	\begin{align}
		\max_{\alpha} \Big \{ f^a(\alpha;x) : a^\top \alpha \le A, \alpha \in \{0,1\}^{|V|} \Big\},
	\end{align}
	where $d \in \mathbb{R}^{|V|}_+$, $a\in \mathbb{R}^{|V|}_+$, $D \in \mathbb{R}_+$, $A \in \mathbb{R}_+$.
	\label{def:CNG}
\end{definition}
\begin{remark}
	In \cref{def:CNG}, we implicitly assume that the defender maximizes a function $f^d$ of $x$ parameterized in $\alpha$ that represents the network's connectivity (or the operativeness) $G$. Symmetrically, the attacker maximizes a function $f^a$ of $\alpha$ parameterized in $x$ that represents the network disruption on $G$. Furthermore, the weights $a$ and $d$ represent the resources the player spends for selecting (\ie, defending or attacking) the nodes, while $A$ and $D$ represent the players' resources budgets.
	When we say the game is \emph{simultaneous}, we mean that the players chose their strategy (\ie, solve their optimization problems) simultaneously without knowing the other players' strategies. When we say the game has \emph{complete information}, we mean that both players are aware of the optimization problems, \ie, they are both aware\footnote{How realistic is the assumption of awareness is discussed in Section \ref{sec:cng}.} of $f^a$, $f^d$, $G$, $a$, $d$, $A$, and $D$.
\end{remark}

We define any feasible $x$ (resp. $\alpha$) as a strategy for the defender (resp. attacker). We call any tuple $(x,\alpha)$ a strategy profile for the game, and $f^d$ (resp. $f^a$) evaluated at $(x,\alpha)$ the defender's (resp. attacker's) payoff under $(x,\alpha)$.
We consider the \gls{NE} as the standard solution concept.
In plain English, a strategy profile ($\bar x,\bar \alpha$) is a \gls{NE} if no player can unilaterally deviate from their strategy without decreasing their payoff. We formalize this concept in \cref{def:NE}.
\begin{definition}[Nash Equilibrium]
	A profile $(\bar x,\bar \alpha)$ is a \gls{NE} if (i.) $f^d(\bar x; \bar \alpha) \ge f^d(\tilde x; \bar \alpha)$ for any $\tilde x \in \{ x \in \{0,1\}^{|V|}: d^\top x \le D\} $, and (ii.) $f^a( \bar \alpha;\bar x) \ge f^a(\tilde \alpha;\bar x)$ for any $\tilde \alpha \in \{ \alpha \in \{0,1\}^{|V|} : a^\top \alpha \le A\} $.
	\label{def:NE}
\end{definition}
We focus on deterministic, or \emph{pure}, \gls{NE}; namely, we assume the players select strategies that are feasible for their constraints instead of selecting convex combinations of feasible strategies. Pure \gls{NE}s may not exist for finite games such as the \gls{CNG}, namely, for games with a finite number of players and a finite number of strategies; indeed, the problem of determining if an \gls{NE} exists in integer programming games is $\Sigma^p_2$-complete \citep{vaz_existence_2018,Dragotto_2021_ZERORegrets}, \ie, as long as $\mathcal{NP} \neq \Sigma^p_2$ the decision problem cannot be represented as an integer program of polynomial size. We refer to \citet{IPGs_2023_Tutorial} for a review of the topic.
Whenever a \gls{NE} does not exist, we rely on the relaxed concept of approximate \gls{NE} that we formalize in \cref{def:ApproximateNE}.

\begin{definition}[Approximate \gls{NE}]
	Given $\Phi \in \mathbb{R}_+$, a strategy profile $(\bar x,\bar \alpha)$ is a (pure) approximate $\Phi-$\gls{NE}  if (i.) $f^d(\bar x, \bar \alpha) + \Phi \ge f^d(\tilde x, \bar \alpha)$ for any $\tilde x \in \{ x \in \{0,1\}^{|V|}: d^\top x \le D\} $, and (ii.) $f^a(\bar \alpha;\bar x) + \Phi \ge f^a( \tilde \alpha;\bar x)$ for any $\tilde \alpha \in \{ \alpha \in \{0,1\}^{|V|} : a^\top \alpha \le A\} $.
	\label{def:ApproximateNE}
\end{definition}

In a $\Phi-$\gls{NE}, the constant $\Phi \in \mathbb{R}_+$ represents an upper bound on the deviations that the players' payoff can have. Whenever $\Phi=0$, the $\Phi-$\gls{NE} is also an exact \gls{NE}.

\subsection{Prices}
Whenever multiple equilibria exist, their properties (\eg, the players' payoffs under the equilibria) may differ; this is precisely why we aim to select the equilibria exhibiting some desired properties. Specifically, we aim to select the \gls{NE} that maximizes the defender's payoff and the \gls{NE} that maximizes the attacker's payoff.  If different, these two \gls{NE}s will practically provide the best “stable” outcome -- in terms of aggressive or defensive strategy -- for the attacker or the defender, respectively. This information can guide the design of extra layers of security mechanisms and help provide real-time prescriptive strategies to defend critical infrastructure. In this section, we formalize these intuitions with the concepts of \gls{POS} and \gls{POA}.
\begin{definition}[Joint Outcomes Space]
	The joint outcomes space for the \gls{CNG} is the set
	\begin{align}
		\mathcal{J} = \big \{ (x,\alpha) : d^\top x \le D,
		a^\top \alpha  \le A,
		\alpha \in \{0,1\}^{|V|}, x \in \{0,1\}^{|V|} \}.
	\end{align}
	\label{def:JOS}
\end{definition}
\vspace{-1em}
The set $\mathcal{J}$ contains all the \emph{outcomes} of the \gls{CNG}, namely, all the feasible strategy profiles that players can play.
If we optimize a real-valued function $g(x,\alpha):\{0,1\}^{|V|^2} \rightarrow \mathbb{R}$ over $\mathcal{J}$, we obtain the strategy profile $(\bar x,\bar \alpha)$ that maximizes $g$. There is no guarantee that the profile $(\bar x,\bar \alpha)$ is a \gls{NE}.
Nevertheless, comparing $(\bar x,\bar \alpha)$ with the \gls{NE}s maximizing $g$ enables us to practically evaluate the \emph{loss in performance} due to the equilibrium conditions. In other words, it enables us to evaluate how much the \emph{stability} conditions of equilibria degrade the value of $g$. This approach, similar to the worst-case analysis for algorithms, was pioneered by \citet{goos_worst-case_1999}. The authors focused on $g$ being the so-called \emph{social welfare}, namely, the sum of the player's payoff; they introduced the Price of Stability, \ie, the ratio between the social welfare of the best possible outcome and the social welfare of the best-possible \gls{NE}.
With a similar spirit to \citep{goos_worst-case_1999}, we introduce a metric for the \gls{CNG} called the \gls{POS} in \cref{def:PoS}.
For a given $\Phi$, let $\mathcal{N}(\Phi) \subseteq \mathcal{J}$ be the \emph{set of $\Phi$-\gls{NE}} for the \gls{CNG}.

\begin{definition}[Price of Security]
	Given a \gls{CNG} instance and $\Phi \in \mathbb{R}_+$, let $(\hat x, \hat \alpha) = \arg \max_{x,\alpha}\{ f^d(x,\alpha) : (x,\alpha) \in \mathcal{N}(\Phi) \}$, \ie, the \gls{NE} maximizing $f^d(x,\alpha)$; let $(\bar x, \bar \alpha)=\arg \max_{x,\alpha}\{ f^d(x,\alpha) : (x,\alpha) \in \mathcal{J} \}$, \ie, the strategy profile maximizing $f^d(x,\alpha)$. Whenever $|\mathcal{N}(\Phi)|>0$, the \gls{POS} is the ratio $f^d(\bar x, \bar \alpha)/f^d(\hat x, \hat \alpha)$.
	\label{def:PoS}
\end{definition}

The \gls{POS} is lower bounded by $1$ and, from a theoretical perspective, has no upper bound. A \gls{POS} of $1$ suggests that the defender, in the best-possible \gls{NE} that maximizes their payoff, is not diminishing their payoff by defending their resources; hence, it suggests that the defensive strategy is highly efficient. A \gls{POS} strictly greater than $1$ suggests a loss in the defender's payoff caused by the \gls{NE} conditions. In this sense, a larger \gls{POS} indicates that the defender is paying a higher cost for defending their resources. We remark that in this case, $f^d(x,\alpha)$ is a real-valued function in $x$ and $\alpha$ and not a function parametrized in $\alpha$.
Symmetrically to the \gls{POS}, we introduce the \gls{POA} in \cref{def:PoA}.

\begin{definition}[Price of Aggression]
	Given a \gls{CNG} instance and $\Phi\in \mathbb{R}_+$, let $(\hat x, \hat \alpha) = \arg \max_{x,\alpha}\{ f^a(\alpha,x) : (x,\alpha) \in \mathcal{N}(\Phi) \}$, \ie, the \gls{NE} maximizing $f^a(\alpha,x)$; let $(\bar x, \bar \alpha)=\arg \max_{x,\alpha}\{ f^a(\alpha,x) : (x,\alpha) \in \mathcal{J} \}$, \ie, the strategy profile maximizing $f^a(\alpha,x)$. Whenever $|\mathcal{N}(\Phi)|>0$, the \gls{POA} is the ratio $f^a(\bar \alpha,\bar x)/f^a(\hat \alpha,\hat x)$.
	\label{def:PoA}
\end{definition}

From the attacker's perspective, the \gls{POA} mimics the definition of the \gls{POS}. It measures the relative loss in the attacker's payoff caused by the \gls{NE} conditions in the \gls{NE} that maximizes the real-valued function $f^a(\alpha,x)$.

\section{Cloud Critical Node Game }
\label{sec:cng}

This section contextualizes our mathematical formulation for the \gls{CNG} in cloud networks. In \cref{sub:payoffs}, we describe the payoff structure for the game. Specifically, we describe $f^d$ and $f^a$. In \cref{sub:computing}, we recall the \emph{ZERO Regrets} algorithm \citep{Dragotto_2021_ZERORegrets} and describe how to tailor it for our problem.

\paragraph{Contextualization.} A malicious attacker gains access to the network $G$ via undisclosed vulnerabilities and aims to perform an attack to maximize the network disruption on $G$ (\eg, maximizing $f^a$). The defender detects the intrusion, yet, they are uncertain about the attacker's targets and the type of vulnerability the attacker could exploit. Therefore, the defender aims to protect their critical infrastructure. The defender creates a projection of the attacker's capabilities by defining the attacker's optimization model, \eg, by choosing $f^a$, $A$, and $a$ based on the monitoring data that warned about a possible attack and the available information regarding the exploited vulnerabilities. The attacker and the defender have limited resources for their defensive and attacking strategies (namely, $a,A,d,D$). While the attacker aims to maximize a measure of network disruption, the defender aims to preserve the network operations as much as possible by maximizing $f^d$. The defender may activate an extra layer of security measures (\ie, a firewall or denial-of-service attack mitigation) on a subset of critical nodes by degrading the network operations by a given factor.

\subsection{The Payoffs}
\label{sub:payoffs}
Let $p^d_i$ be the parameter representing the criticality of node $i \in V$ for the defender. Symmetrically, let $p^a_i$ be the parameter representing the projected criticality of node $i$ according to the attacker, \eg, according to the vulnerabilities the attacker could exploit and the importance of $i$. We will formulate the payoffs for the two players in terms of $p^d_i$ and $p^a_i$ according to whether node $i$ is protected ($x_i=1$) or attacked ($\alpha_i=1$), respectively. Let $\delta$, $\eta$, $\epsilon$ and $\gamma$ be real-valued scalar parameters in $[0,1]$ so that $\delta < \eta < \epsilon$. The payoff contributions follow the following scheme:

\begin{enumerate}
	\item \textbf{Normal operations.} If $x_i=0$ and $\alpha_i=0$, the defender gets a full payoff of $p^d_i$ as no attack is ongoing on $i$. However, the attacker pays an opportunity cost $\gamma p^a_i$ for not attacking $i$.

	\item \textbf{Successful attack.} If $x_i=0$ and $\alpha_i=1$, the attacker gets a full payoff of $p^a_i$, as the attacker successfully attacked node $i$. Therefore, the defender's operations on node $i$ worsen from $p^d_i$ to $\delta p^d_i$.

	\item \textbf{Mitigated attack.} If $x_i=1$ and $\alpha_i=1$, the defender's selection of node $i$ mitigates the attack. Therefore, the defender's operations are degraded from $p^d_i$ to $\eta p^d_i$. Symmetrically, the attacker receives a payoff of $(1-\eta) p^a_i$.

	\item \textbf{Mitigation without attack.} If $x_i=1$ and $\alpha_i=0$, the defender protects node $i$ without the attacker selecting node $i$. Therefore, the defender's operations are degraded from $p^d_i$ to $\epsilon p^d_i$. Symmetrically, the attacker receives a payoff of $0$.
\end{enumerate}

We summarize the contribution of each node $i$ to the players' payoff in \cref{tab:payoff}.
\begin{table}[!ht]
	\centering
	\begin{tabular}{l@{\hspace{2em}}|c@{\hspace{2em}}c@{\hspace{2em}}}
		             & $\bm{\alpha_i=0}$                                       & $\bm{\alpha_i=1}$                                            \\
		\midrule
		$\bm{x_i=0}$ & $ {\color{blue}p^d_i}$ $|$ ${\color{red}-\gamma p^a_i}$ & $ {\color{blue}\delta p^d_i}$ $|$ ${\color{red} p^a_i}$      \\
		$\bm{x_i=1}$ & $ {\color{blue}\epsilon p^d_i}$ $|$ ${\color{red}0}$    & $ {\color{blue}\eta p^d_i}$ $|$ ${\color{red}(1-\eta)p^a_i}$
	\end{tabular}
	\caption{Description of the contributions to the player's payoff with respect to the strategies associated with each node $i$. In blue (resp. red), are the defender's (resp. attacker's) contributions to their payoff.\label{tab:payoff}}
\end{table}

All considered, the defender's payoff is
\begin{align}
	f^d(x,\alpha)=\sum_{i \in V} \Big (p^d_i\big ( (1-x_i)(1-\alpha_i) + \eta x_i\alpha_i  +\epsilon x_i(1-\alpha_i) + \delta(1-x_i)\alpha_i \big) \Big ),
\end{align}
while the attacker's payoff is
\begin{align}
	f^a(\alpha,x)=\sum_{i \in V} \Big (p^a_i\big ( -\gamma(1-x_i)(1-\alpha_i) + (1-x_i)\alpha_i  +(1-\eta)x_i\alpha_i \big) \Big ).
\end{align}

\begin{example}
	Consider a cloud network given by $G=(V,E)$ with $|V|=6$, and $|E|=7$ as in \cref{fig:example}. The cloud operator manages nodes $1,\dots,5$ while node $6$ virtually represents the external network (\eg, the Internet). An attacker penetrates the network $G$ and observes the traffic exchanged among the nodes. Once the attack has been detected, the network operator solves a \gls{CNG} to determine the best defensive strategy. Based on the information received from the monitoring systems, the network operator formulates a \gls{CNG} with the parameters of \cref{tab:example}.
	\begin{figure}[!htpb]
		\centering
		\includegraphics[width=0.6\linewidth]{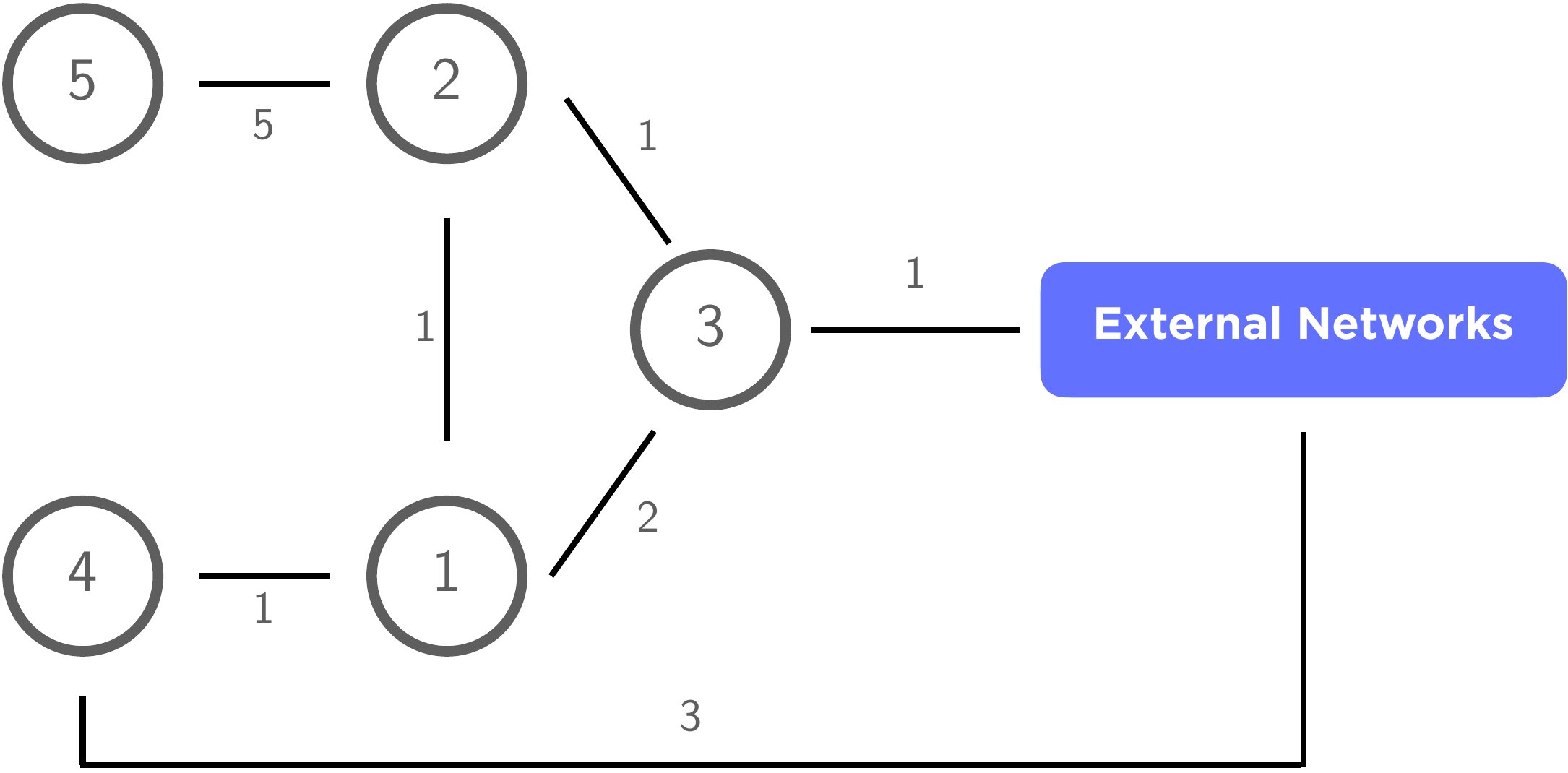}
		\caption{The network of \cref{ex:example}. The weights on the edges represent the amount of traffic exchanged between the nodes.}
		\label{fig:example}
	\end{figure}

	\begin{table}[!ht]
		\centering
		\begin{tabular}{lrrrrr}
			\toprule
			                 & \textbf{Node 1} & \textbf{Node 2} & \textbf{Node 3} & \textbf{Node 4} & \textbf{Node 5} \\
			\midrule
			\textbf{Traffic} & 4               & 7               & 4               & 4               & 5               \\
			$\bm{p^d_i}$     & 9               & 2               & 30              & 3               & 8               \\
			$\bm{p^a_i}$     & 6               & 10.5            & 18              & 6               & 7.5             \\
			$\bm{d}$         & 3               & 6               & 8               & 7               & 7               \\
			$\bm{a}$         & 6               & 4               & 7               & 9               & 1               \\
			\bottomrule
		\end{tabular}
		\qquad
		\begin{tabular}{lr}
			\toprule
			\textbf{Parameter} & \textbf{Value} \\
			\midrule
			$A$                & $25.50$        \\
			$D$                & $40.00$        \\
			$\delta$           & $0.06$         \\
			$\eta$             & $0.40$         \\
			$\epsilon$         & $1.00$         \\
			\bottomrule
		\end{tabular}
		\caption{Parameters associated with the instance of \cref{ex:example}.}
		\label{tab:example}
	\end{table}

	The associated \gls{CNG} admits the two exact (\ie, $\Phi=0$) \gls{NE}s of \cref{tab:exampleNE}. Specifically, the defender achieves, in both cases, a payoff $f^d(\bar{x},\bar{\alpha})$ of $29.2$. However, the attacker has a payoff $f^a(\bar{\alpha},\bar{x})$ of $13.74$ and $12.18$ in \emph{equilibrium $1$} and $2$, respectively. The \gls{POS} for the instance is $1.78$ for both equilibria (\eg, the ratio between $\sum_i p^d_i = 52$ and $29.2$), and it intuitively explains that the defender degrades their payoff $1.78$ times to defend the network with the Nash equilibria. Symmetrically, the \gls{POA} for the instance is $1.86$ and $2.09$ in \emph{equilibrium} $1$ and $2$, respectively; this price indicates that the attacker's payoff decreases at most $2.09$ times in \emph{equilibrium $2$} compared to the most successful attack (\eg, when the defender does not protect any node).
	Although the two equilibria are equally advantageous for the defender, the \gls{POA} provides an insight into the best equilibrium. Indeed, from a practical perspective, the network operator should select \emph{equilibrium $2$}, as it damages the attacker the most.

	\begin{table}[!ht]
		\centering
		\begin{tabular}{cccll}
			\toprule
			              & $\bm{\bar{x}}$ & $\bm{\bar{\alpha}}$ & \gls{POS} & \gls{POA} \\
			\midrule
			Equilibrium 1 & $[1,1,1,0,1]$  & $[0,0,1,0,1]$       & $1.78$    & $1.86$    \\
			Equilibrium 2 & $[0,1,1,0,1]$  & $[0,0,1,0,1]$       & $1.78$    & $2.09$    \\
			\bottomrule
		\end{tabular}
		\caption{The \gls{NE} associated with the instance of \cref{ex:example}.}
		\label{tab:exampleNE}
	\end{table}

	\label{ex:example}
\end{example}

\begin{remark}[Determining the criticality of the nodes]
	In general, there are no standard and shared guidelines on choosing the criticality of each node in an “exact” fashion, as this depends on the characteristics of the network and the subjective assessment of the potential threats associated with each node.
	For instance, nodes that handle sensitive data and authentication mechanisms are more critical due to their impact on data integrity and security. Similarly, a node responsible for authentication is crucial in ensuring the application's security. If this component is compromised or fails, it can lead to unauthorized access and potential data breaches. The criticality of nodes can vary depending on the specific context and requirements of the cloud-native applications. Factors such as business objectives, customer expectations, and system and network dependencies can all contribute to determining the criticality of each node.
\end{remark}

\subsection{Computing Equilibria}
\label{sub:computing}

This section describes the algorithm we employ to compute the equilibria for the \gls{CNG}.
We tailor the cutting-plane algorithm \emph{ZERO Regrets} from \citet{Dragotto_2021_ZERORegrets} to compute \gls{NE}s and $\Phi-$\gls{NE} in the \gls{CNG}. The algorithm receives as an input the \gls{CNG} instance and a function $f(x,\alpha)$ and returns the \gls{NE} maximizing $f$.
We remark that \emph{ZERO Regrets} can compute, thanks to its theoretical guarantees, the equilibrium maximizing $f$; therefore, our algorithm inherits \emph{ZERO Regrets}'s properties and can provably return the equilibrium that maximizes $f(x,\alpha)$, or the approximate equilibrium that maximizes $f(x,\alpha)$ for a given approximation constant $\Phi_{UB}$. Our algorithm will return an equilibrium independently of any network parameter or network topology and can even handle arbitrary linear constraints besides the knapsack constraint.
\SetKwBlock{Repeat}{repeat}{}
\begin{algorithm}[!ht]
	\DontPrintSemicolon
	\caption{ZERO Regrets for the \gls{CNG}\label{Alg:ZERO}}
	\KwData{A \gls{CNG} instance and a function $f(x,\alpha):\{0,1\}^{|V|^2}\rightarrow\mathbb{R}$.}
	\KwResult{ Either: (i.) the \gls{NE} $\bar{x},\bar{\alpha}$ maximizing $f(x,\alpha)$, or (ii.)  a $\Phi-$\gls{NE} $\bar{x},\bar{\alpha}$ maximizing $f(x,\alpha)$. }
	$\Omega =\{ 0 \le 1\}$, $\Phi_{UB}=0$, and $\mathcal{Q}=\max_{x, \alpha, \Phi} \{ f(x,\alpha) : (x,\alpha) \in \mathcal{J}, (x,\alpha,\Phi) \in \Omega, \Phi \le \Phi_{UB} \}$\; \label{Alg:ZERO:Init}
	\Repeat{
		\lIf{$\mathcal{Q}$ is infeasible}
		{$\Phi_{UB}=\Phi_{UB}+1$ \label{Alg:ZERO:Increment}}
		\Else{
			$(\bar x, \bar \alpha, \bar \Phi) \in \arg\max \mathcal{Q}$; \;
			$\tilde x \in \arg\max_x\{f^d(x,\bar{\alpha}) : d^\top x \le D\}$\; \label{Alg:ZERO:BRd}
			$\tilde \alpha \in \arg\max_\alpha\{f^a(\alpha,\bar{x}) : a^\top \alpha \le A\}$\;\label{Alg:ZERO:BRa}
			\If{$f^d(\bar{x},\bar{\alpha}) + \Phi_{UB} \le f^d(\tilde{x},\bar{\alpha})$ \label{Alg:ZERO:CheckD}}
			{
				add $f^d(\tilde{x},\alpha) \le f^d(x,\alpha) + \Phi$ to $\Omega$ \label{Alg:ZERO:CutD}
			}
			\ElseIf{$f^a(\bar{\alpha},\bar{x}) + \Phi_{UB} \le f^a(\tilde{\alpha},\bar{x})$}
			{
				add $f^a(\tilde \alpha,x) \le f^a(\alpha,x)+ \Phi$ to $\Omega$ \label{Alg:ZERO:CutA}
			}
			\Else{
				\KwRet{the $\bar{\Phi}$-\gls{NE} $\bar{x},\bar{\alpha}$} \label{Alg:ZERO:ReturnPNE}
			}
		}
	}
\end{algorithm}

\paragraph{The Algorithm.} We formalize the algorithm in \cref{Alg:ZERO}. In \cref{Alg:ZERO:Init}, the algorithm initializes (i.) an empty cutting plane pool $\Omega$, and (ii.) $\Phi$ to $0$, and (iii.) a program $\mathcal{Q}$ that optimizes the input function $f(x,\alpha)$ over the joint outcome space and $\Omega$. Without loss of generality, we assume that $\mathcal{Q}$ is feasible and bounded at the first iteration (see \cite{Dragotto_2021_ZERORegrets} for details).
Let $\bar x, \bar \alpha$ be the maximizers of $\mathcal{Q}$. The task is to determine whether or not $\bar x, \bar \alpha$ is an \gls{NE}. This is equivalent to checking whether the defender (resp. the attacker) can \emph{unilaterally} and \emph{profitably} deviate to another strategy $\tilde x$ (resp. $\tilde \alpha$). The algorithm solves the defender's (resp. attacker's) optimization problem in \cref{Alg:ZERO:BRd}  (resp. \cref{Alg:ZERO:BRa} for the attacker), by letting $x$ (resp. $\alpha$) be a variable while fixing $\alpha$ to $\bar \alpha$ (resp. $x$ to $\bar x$). If the defender payoff $f^d$ (resp. attacker payoff $f^a$) under the profile $\tilde x, \bar \alpha$ (resp. $\bar x$, $\tilde \alpha$) is better than the one under $\bar x, \bar \alpha$, then $\tilde x$ (resp. $\tilde \alpha$) is a deviation; by definition, $\bar x, \bar \alpha$ cannot be a \gls{NE}. Therefore, the algorithm cuts off $\bar x, \bar \alpha$ from $\mathcal{Q}$ via the cutting plane in \cref{Alg:ZERO:CutD} (resp. \cref{Alg:ZERO:CutA}) with a so-called \emph{equilibrium inequality} \citep{Dragotto_2021_ZERORegrets}, \ie, an inequality that does not cut off any \gls{NE}.
If no deviation exists for the attacker and the defender, then the algorithm returns a \gls{NE} in \cref{Alg:ZERO:ReturnPNE}.
We generalize the previous reasoning by allowing any profitable deviation to be incremented by at most $\bar{\Phi}$ to enable the computation of $\Phi-$\gls{NE}. Whenever $\mathcal{Q}$ becomes infeasible (\cref{Alg:ZERO:Increment}), then no $\Phi$-\gls{NE} can exist with $\Phi \le \Phi_{UB}$; therefore, we heuristically increment $\Phi_{UB}$ by one unit.

\begin{remark}
	The network operators run \cref{Alg:ZERO} with a time limit and expect the algorithm to produce a feasible solution. In order to guarantee that the algorithm produces a feasible solution, at each iteration, we store the best incumbent $\Phi$-\gls{NE} found so far, that is, the solution $\bar x, \bar \alpha$ with the smallest $\bar{\Phi}$. If the algorithm hits a time limit, we return the best incumbent $\Phi$-\gls{NE} found.
\end{remark}

\section{Computational Experiments}
\label{sec:tests}
We perform the computational experiments on two instance sets. Specifically, we test our algorithm on synthetic instances and instances derived from a real-world cloud network. We run our experiments on $8$ CPUs and $64$ GB of RAM, with \emph{Gurobi 10} as the optimization solver.
The full tables of results and instances are available at \url{https://github.com/ds4dm/CNG-Instances}.
\subsection{Synthetic Instances}
We generate a series of cloud networks with $10,25,50,75,100,150$ and $300$ nodes.
We draw the parameters $\delta, \eta, \epsilon, \gamma, A$, and $D$ from a series of realistic distributions mimicking the dynamics of several attack loads.  For instance, low-load attacks represent activities related to smaller attacks, while high-load attacks represent attacks of major impact (\eg, distributed attacks, state-sponsored attacks). We summarize the choices for such parameters in \cref{tab:params}.

\begin{table}[ht]
	\centering
	\begin{tabular}{@{}c@{\hspace{1em}}r@{\hspace{2em}}l@{\hspace{1em}}}
		\toprule
		\textbf{Parameter} & \textbf{Values}                                        & \textbf{Notes}                              \\
		\midrule
		$\gamma$           & $0.00, 0.10$                                           & Attacker's opportunity cost factor          \\
		$\eta$             & $0.60, 0.80$                                           & Defender's mitigated-attack factor          \\
		$\epsilon$         & $1.25\eta$                                             & Defender's mitigation-without-attack factor \\
		$\delta$           & $0.80\eta$                                             & Attacker's successful-attack factor         \\

		$D$                & $0.30\sum_{i} d_i, 0.75\sum_{i} d_i$                   & Defender's budget                           \\
		$A$                & $0.03\sum_{i} a_i, 0.10\sum_{i}a_i , 0.30\sum_{i} a_i$ & Attacker's budget                           \\
		\bottomrule
	\end{tabular}
	\label{tab:parameters}
	\caption{Description of the \gls{CNG} parameters for the synthetic instances. \label{tab:params}}
\end{table}

The motivation behind the values of the parameters is mainly empirical. The attacker's opportunity cost $\gamma$ is generally considered to be either $0$ or $0.1$ as attacking may expose the exact dynamics of the attack (\eg, the type of vulnerabilities exploited by the attacker) while not attacking may give time to the network operator to investigate on the attack. The magnitude of the mitigated attack factor $\eta$ depends on the type of defensive resources deployed by the network operator; for instance, some firewall filtering rules may significantly slow down the overall network operations, whereas some may have milder effects on the network performance. Therefore, we select a low (\ie, $0.60$) and a high (\ie, $0.80$) value of $\eta$. Since no attack is ongoing, the mitigation-without-attack factor $\epsilon$ is based on a $\eta$ plus an extra benefit (\ie, $\epsilon = 1.25 \eta$). The successful-attack factor $\delta$ is based on $\eta$ minus an extra cost due to the unmitigated attack (\ie, $\delta = 0.8 \eta$). The budgets $A$ and $D$ are strictly instance-dependent since they model the players' ability to select nodes in the network. In our instances, we consider large-scale attacks where the defender generally has resources to protect either $30\%$ or $75\%$ of the network, as of contractual agreements with the customers regarding the degradation of services (\ie, the so-called \emph{service level agreements}). As of these agreements, the defender may be contractually obliged to guarantee a minimal level of service to its customer. In contrast, the attacker can attack the $3\%,10\%$ or $30\%$ of the nodes according to their weights $a$.
We generate $a,d$ as random integer vectors with entries in the range $[1,25]$. We generate $p^a$ and $p^d$ by starting from the same random integer vector with entries in the range $[1,25]$ and adding, for $p^a$ and $p^d$ separately, another random integer vector with the same characteristics.

\begin{table}[!ht]
	\centering
	\resizebox{\textwidth}{!}{%
		\begin{tabular}{c@{\hspace{3em}}r@{\hspace{2em}}r@{\hspace{2em}}r@{\hspace{2em}}r@{\hspace{2em}}r@{\hspace{2em}}r@{\hspace{2em}}r@{\hspace{2em}}r@{\hspace{2em}}}
			\toprule
			$\bm{|V|}$ & \gls{POS} & \gls{POS} \textbf{Range} & \gls{POA} & \gls{POA} \textbf{Range} & $\bm{\Phi}$ & $\bm{f^d}$ & $\bm{f^a}$ & \textbf{Time (s)} \\
			\midrule
			10         & 1.10      & {[}1.00, 1.39{]}         & 2.59      & {[}1.38, 5.00{]}         & 15.55       & 1731.45    & 99.92      & 27.65             \\
			25         & 1.09      & {[}1.01, 1.34{]}         & 2.72      & {[}1.43, 5.00{]}         & 4.53        & 729.41     & 35.38      & 12.24             \\
			50         & 1.11      & {[}1.00, 1.35{]}         & 2.17      & {[}1.56, 2.98{]}         & 19.63       & 1592.66    & 90.34      & 34.08             \\
			75         & 1.11      & {[}1.00, 1.36{]}         & 2.40      & {[}1.07, 3.24{]}         & 22.57       & 2211.23    & 129.07     & 48.87             \\
			100        & 1.12      & {[}1.00, 1.30{]}         & 4.27      & {[}1.60, 7.72{]}         & 41.60       & 3152.46    & 157.07     & 65.57             \\
			150        & 1.17      & {[}1.01, 1.39{]}         & 5.02      & {[}1.51, 10.16{]}        & 75.71       & 6842.50    & 272.88     & 86.69             \\
			300        & 1.10      & {[}1.01, 1.26{]}         & 4.09      & {[}1.42, 7.60{]}         & 174.82      & 8592.27    & 479.91     & 100.02            \\
			\bottomrule
		\end{tabular}
	}
	\caption{Average computational results for the synthetic instances. The results are aggregated for each value of cardinality of $V$.}
	\label{tab:results_syn}
\end{table}

\paragraph{Results.} In \cref{tab:results_syn}, we present the computational results for the synthetic instances. Each row presents some metrics averaged over all instances with the same number of nodes $|V|$. For each combination of parameters described in \cref{tab:parameters}, we execute \cref{Alg:ZERO} twice by setting $f(x,\alpha)=f^d(x,\alpha)$ and $f(x,\alpha)=f^a(\alpha,x)$, \ie, we execute the algorithm twice to compute the \gls{POA} and \gls{POS}.
In column order, we report
\begin{enumerate*}
	\item the number of nodes $|V|$,
	\item the average \gls{POS}
	\item the range of \gls{POS}, \ie, the minimum and maximum \gls{POS} achieved in the instances,
	\item the average \gls{POA},
	\item the range of \gls{POA},
	\item the average value of $\Phi$ for $\Phi$-\gls{NE},
	\item the average defender's payoff $f^d$,
	\item the average attacker's payoff $f^a$, and
	\item the average time (seconds) required to compute the \gls{NE}
\end{enumerate*}.
In all the instance sets but $300$, we compute a \gls{NE} within the time limit of $100$ seconds. In the set with $|V|=300$, the algorithm cannot compute an equilibrium within the time limit and returns a strategy that, on average, is a $\Phi$-\gls{NE} with $\Phi=174.82$, whereas $f^d$ and $f^a$ average $8592.27$ and $479.91$, respectively. This indicates that, on average, $\Phi$ is $2\%$ of the defender's payoff, and thus the $\Phi$-\gls{NE} has a small relative approximation ratio. The equilibria computed in all the instances tend to exhibit similar behavior in terms of \gls{POS}, whereas the \gls{POA} tends to grow in the sets $100,150$ and $300$. This growth in the \gls{POA} is also partially reflected in the broader range in the instances $100, 150$, and $300$. We remark that lower \gls{POS} indicate that the defender can commit to equilibria strategies that are close, in terms of efficiency, to the best possible strategy in the joint outcome space. In other words, equilibria with low \gls{POS} are practically-efficient defensive strategies.

\subsection{Real-world Cloud Network}
We instantiate the \gls{CNG} on a real-world anonymized dataset from Ericsson \citep{boukhtouta_cloud_2022} containing traffic snapshots from a cloud-native application network. Specifically, we consider $10$ snapshots of the same network collected at different times (\ie, each $1$ minute), and we perform our tests on each of the snapshots. The network is a cluster node running a set of cloud services (\eg, containers registry, logging, database mediator, key management, load balancer, DNS, LDAP, SNMP, and SCTP services). In addition, there are $3$ management hosts (\eg, Kubernetes K8S), and $16$ servers providing customer resources. This cloud network powers $|V|=656$ separate resources (\eg, storage, virtual servers). We refer the reader to \citep{boukhtouta_cloud_2022} for additional information concerning the data source and the underlying cloud network.
We assign the parameters $p^d_i$ and $p^a_i$ by summing up the total traffic transiting through node $i$. Similarly, we derive the nodes' weights $d_i$ and $a_i$ as $\log_2(p^d_i)$. We modify the defender's $p^d_i$ and $d_i$ to reflect the node's nature, \ie, we increase them if they are associated with critical infrastructures, such as management nodes (\eg, K8S master node, DNS service discovery or service mesh routing) and vital nodes for the cloud network. As a result of this latter altering procedure, we remark that the final weights' magnitude might not be correlated with the total traffic of management nodes. Similarly, we modify the attacker's profits $p^a_i$ and weights $a_i$ to partially reflect the attacker's knowledge of the nature of the node, \ie, an increased cost on $i$ may signal that the attacker is more confident that the security mechanisms at node $i$ is difficult to penetrate.

\begin{table}[!ht]
	\centering
	\resizebox{\textwidth}{!}{%
		\begin{tabular}{llllllrrrrrrrrrr}
			\toprule
			\textbf{A\%} & \textbf{D\%} & $\bm{\eta}$ & $\bm{\epsilon}$ & $\bm{\gamma}$ & $\bm{\delta}$ & \gls{POS} & \gls{POS} \textbf{Range} & \gls{POA} & \gls{POA} \textbf{Range} & $\bm{\Phi}$ & $\bm{f^d}$ & $\bm{f^a}$ & \textbf{Time (s)} \\
			\midrule
			0.03         & 0.30         & 0.60        & 0.75            & 0.00          & 0.48          & 1.07      & {[}1.05, 1.10{]}         & 1.04      & {[}1.01, 1.09{]}         & 92.10       & 8366.01    & 3677.03    & 72.70             \\
			0.10         & 0.30         & 0.60        & 0.75            & 0.00          & 0.48          & 1.12      & {[}1.11, 1.13{]}         & 1.01      & {[}1.00, 1.04{]}         & 241.33      & 7948.18    & 9429.20    & 180.00            \\
			0.03         & 0.75         & 0.60        & 0.75            & 0.00          & 0.48          & 1.07      & {[}1.05, 1.10{]}         & 1.04      & {[}1.01, 1.09{]}         & 92.10       & 8365.98    & 3677.03    & 74.73             \\
			0.10         & 0.75         & 0.60        & 0.75            & 0.00          & 0.48          & 1.12      & {[}1.11, 1.13{]}         & 1.01      & {[}1.00, 1.04{]}         & 241.33      & 7948.18    & 9429.20    & 180.00            \\
			0.30         & 0.75         & 0.60        & 0.75            & 0.00          & 0.48          & 1.24      & {[}1.21, 1.26{]}         & 1.13      & {[}1.01, 1.26{]}         & 394.74      & 7117.08    & 16106.87   & 180.00            \\
			0.03         & 0.30         & 0.80        & 1.00            & 0.00          & 0.64          & 1.01      & {[}1.01, 1.04{]}         & 3.39      & {[}2.33, 4.52{]}         & 57.34       & 8767.94    & 971.60     & 180.00            \\
			0.10         & 0.30         & 0.80        & 1.00            & 0.00          & 0.64          & 1.08      & {[}1.07, 1.09{]}         & 1.48      & {[}1.00, 3.41{]}         & 284.97      & 8278.47    & 8769.75    & 180.00            \\
			0.03         & 0.75         & 0.80        & 1.00            & 0.00          & 0.64          & 1.01      & {[}1.01, 1.03{]}         & 3.35      & {[}2.44, 5.00{]}         & 51.80       & 8774.46    & 996.21     & 180.03            \\
			0.10         & 0.75         & 0.80        & 1.00            & 0.00          & 0.64          & 1.08      & {[}1.07, 1.09{]}         & 1.01      & {[}1.00, 1.04{]}         & 307.24      & 8248.47    & 9429.20    & 180.00            \\
			0.30         & 0.75         & 0.80        & 1.00            & 0.00          & 0.64          & 1.16      & {[}1.14, 1.17{]}         & 1.02      & {[}1.01, 1.03{]}         & 541.44      & 7707.95    & 16893.80   & 180.00            \\
			0.03         & 0.30         & 0.60        & 0.75            & 0.10          & 0.48          & 1.07      & {[}1.05, 1.10{]}         & 3.74      & {[}2.72, 6.40{]}         & 92.83       & 8342.57    & 715.32     & 89.77             \\
			0.10         & 0.30         & 0.60        & 0.75            & 0.10          & 0.48          & 1.20      & {[}1.15, 1.23{]}         & 1.74      & {[}1.23, 2.03{]}         & 676.05      & 7612.48    & 6058.38    & 180.00            \\
			0.03         & 0.75         & 0.60        & 0.75            & 0.10          & 0.48          & 1.07      & {[}1.05, 1.10{]}         & 5.51      & {[}3.82, 10.10{]}        & 92.87       & 8347.80    & 713.36     & 91.69             \\
			0.10         & 0.75         & 0.60        & 0.75            & 0.10          & 0.48          & 1.14      & {[}1.11, 1.15{]}         & 1.89      & {[}1.36, 2.19{]}         & 730.06      & 7806.45    & 4525.63    & 180.00            \\
			0.30         & 0.75         & 0.60        & 0.75            & 0.10          & 0.48          & 1.27      & {[}1.19, 1.41{]}         & 1.32      & {[}1.15, 1.42{]}         & 662.80      & 7014.45    & 12067.67   & 180.00            \\
			0.03         & 0.30         & 0.80        & 1.00            & 0.10          & 0.64          & 1.02      & {[}1.01, 1.05{]}         & 3.88      & {[}1.01, 9.48{]}         & 85.86       & 8693.34    & 910.89     & 180.00            \\
			0.10         & 0.30         & 0.80        & 1.00            & 0.10          & 0.64          & 1.08      & {[}1.07, 1.09{]}         & 1.18      & {[}1.00, 2.17{]}         & 289.00      & 8252.51    & 7955.36    & 180.00            \\
			0.03         & 0.75         & 0.80        & 1.00            & 0.10          & 0.64          & 1.01      & {[}1.01, 1.03{]}         & 4.27      & {[}2.40, 5.25{]}         & 60.03       & 8753.73    & 827.67     & 180.00            \\
			0.10         & 0.75         & 0.80        & 1.00            & 0.10          & 0.64          & 1.08      & {[}1.07, 1.09{]}         & 1.91      & {[}1.02, 3.47{]}         & 249.37      & 8308.51    & 7584.37    & 180.00            \\
			0.30         & 0.75         & 0.80        & 1.00            & 0.10          & 0.64          & 1.16      & {[}1.14, 1.17{]}         & 1.02      & {[}1.01, 1.03{]}         & 541.43      & 7705.90    & 16891.54   & 180.00            \\
			\bottomrule
		\end{tabular}
	}
	\caption{Average computational results for the real-world cloud network.}
	\label{tab:results_real}
\end{table}

\paragraph{Results.} In \cref{tab:results_real}, we report the aggregated computational results conforming to the previous description. We aggregate the instances by their parameters $\eta, \epsilon, \gamma, \delta$. For each instance, we run \cref{Alg:ZERO} twice over all the combinations of parameters %
to determine the \gls{POA} and the \gls{POS}. We employ the same notation for the columns, except for the columns $A\%$ and $D\%$; in these columns, we report the weighted percentage of nodes the attacker or defender can select. For instance, $A\%=0.3$ indicates the attacker's budget is $0.3\sum_i a_i$. Since the instances are generally larger than the synthetic ones, the algorithm almost always hits the time limit of $180$ seconds. However, it also always finds a feasible $\Phi$-\gls{NE} that ranges from $1\%$ to $9\%$ of the defender's payoff, \ie, equilibria that are close to being exact. The \gls{POS} varies according to the type of attack (\ie, with the parameters), and generally degrades the defender's performance from $1\%$ for small-scale attacks to $27\%$ in large-scale attacks, \eg, when $A\% = 0.3$.

\section{Conclusions}\label{sec:conclusion}
In this work, we introduced the \gls{CNG}, a game-theoretic model to assess the cyber-security risk of cloud networks and inform security experts on the optimal security strategies.
Our approach combines game theory and combinatorial optimization to provide a practical framework to guide security experts in assessing the security posture of the cloud network and dynamically adapting the level of cyber-protection deployed on the network.
Methodologically, we formalized the \gls{CNG} as a simultaneous and non-cooperative attacker-defender game where each player solves a combinatorial optimization problem parametrized in the variables of the other player. Our model embeds the uncertainty of the attack dynamics by modeling the range of possible attacks and requires, in this sense, almost no probabilistic assumption.
In practice, we provided a computational analysis of a real-world cloud network and synthetic instances and demonstrated the efficacy of our approach. Importantly, our approach always produces a feasible recommendation for the network operators and scales up to real-sized cloud networks. Although our \gls{CNG} only includes a knapsack constraint, it could support more sophisticated constraints modeling requirements on the defender's and attacker's operations.
Finally, our approach does not assume any specific structure on the attacks (if not a resource constraint). To alleviate the conservativism of such a model, our approach can be augmented by specifying the dynamics (or types) of attacks allowed by the attacker. For instance, the MITRE ATT\&CK framework \citep{mitre-attack}, a collection of attack tactics and techniques based on real-world observations, can provide useful constraints on the attacker dynamics. For example, we could require that if one given node is attacked, other nodes with the same “properties” (\eg, the same software) are also attacked. As long as we can represent these dynamics with linear constraints, our algorithm can still compute the equilibrium of the resulting game.

\section*{Acknowledgements}
We are grateful to Margarida Carvalho for the insightful feedback on this work. We are thankful to the two anonymous referees for their feedback and suggestions.
The third author acknowledges the collaboration and contribution of IASI, CNR, Rome.
\bibliography{Biblio}
\label{sec:bib}

\end{document}